\documentclass{article}

\title{Some results about the projective normality of abelian varieties.}
\author{Luis Fuentes Garc\'{\i}a}
\date{14/06/2004}

\newtheorem{teo}{Theorem}[section]

\newtheorem{cor}[teo]{Corollary}
\newtheorem{lemma}[teo]{Lemma}

\newtheorem{conj}[teo]{Conjecture}

\def\e#1{exp({{#1}})}

\def\g2{\pi}

\def\ZZ{\leavevmode\hbox{$\rm Z$}}

\def\RR{\leavevmode\hbox{$\rm I\!R$}}
\def\CC{\leavevmode\hbox{$\rm I\!\!\!C$}}

\def\qed{\hspace{\fill}$\rule{2mm}{2mm}$}
\newcommand\lrw{\longrightarrow}

\begin{document}

\maketitle

\vspace{0.1cm}

{\bf Abstract:} We reduce the problem of the projective normality of polarized abelian varieties to check the rank of very explicit matrices. This allow us to prove some results on normal generation of primitive line bundles on abelian threefolds and fourfolds. We also give two situations where the projective normality always fails. Finally we make some conjecture. 

\vspace{0.1cm}

\section{Introduction}\label{intro}

Let $(A,L)$ be a polarized abelian variety of dimension $g$. We are interested in studying the projective normality of the embedding given by the linear system $|L|$. The problem is completely solved when $g=1$ or $g=2$ (see \cite{bila}, \cite {la}, \cite{iy1}, \cite{fu}). For higher dimension, most of the existent results are refered to powers of ample line bundle:

\begin{enumerate}

\item If $L=M^n$ with $M$ ample and $n\geq 3$ then $L$ is normally generated (see \cite{bila}).

\item If $L=M^2$ with $M$ of characteristic $c$, then $L$ is normally generated if and only if no point of $t^*_{\bar{c}}K(L)$ is a base point of $M$ (see \cite{bila}). 

\item If $L$ is of type $(2,d_2,\ldots,d_g)$ and $d_g\leq 4$  then $L$ is never normally generated; if $(A,L)$ is generic of type $(2,d_2,\ldots,d_g)$ and  $d_g>4$ then $L$ is normally generated (see \cite{ru}).

\end{enumerate}

On the other hand, the new  techniques of  Mukai-Fourier transform  and M-regularity have been applied by G.Pareschi an M.Popa  to attack these problems (see  \cite{papo1}, \cite{papo2}). However,  it  seems that they only  work well with powers of ample line bundles.

Thus, if  $L$ is primitive, very little is known. In \cite{iy}, J.Iyer proves that if $A$ is simple and $h^0(L)>2^g g!$ then $L$ is normally generated. First she prove that it suffices to study the $2$-normality. Then, she takes a suitable isogeny $A\lrw B$ such that $L$ descends to a principal polarization $M$ on $B$. In this way, the problem is reduced to check the surjectivity of some multiplication maps between translations of the line bundle $M$.

Here, by using the canonical theta functions, we develop this idea. The main theorem reduces the question of $2$-normality to check the rank of some matrices. Their elements are canonical theta functions evaluated in $0$. The most important fact is that we can check the projective normality of any line bundle making a explicit computation. We use this Theorem in two ways.

Firstly, we use the main Theorem to show two situations where the projective normality always fails. If $(A,L)$ is a $g$-dimensional polarized abelian variety of type $(d_1,\ldots,d_g)$ with some $d_i=2$, we prove:

\begin{enumerate}

\item If $d_j\leq 4$, $1\leq j\leq g$  then $L$ is not normally generated. This extends the result of E.Rubei in \cite{ru}.

\item If $h^0(L)=2^{g+1}$, then $L$ is not normally generated. \qed

\end{enumerate}

Secondly, by using a computer we exhibit explicit examples of projectively normal abelian varieties.   Note, that a necessary condition for the projective normality is $h^0(L)\geq 2^{g+1}-1$. Moreover, by the result of J.Iyer we only have to check the cases where $h^0(L)\leq 2^g g!$. Because the projective normality is an open condition, these examples allow us to prove the normal generation in the generic case for some fixed types. For abelian threefolds we prove the following:

\begin{enumerate} 

\item If $(A,L)$ is generic of type $(d_1,d_2,d_3)$ with $h^0(L)\geq 15$ and
$
(d_1,d_2,d_3)\not\in \{(1,2,8),(1,3,6),(2,2,4),(2,4,4)\}
$
then $L$ is  normally generated.

\item If $L$ is of type $(1,2,8)$, $(2,2,4)$ or $(2,4,4)$ then it is never normally generated.

\end{enumerate}

Note, that the open case is when $L$ is of type $(1,3,6)$.

For abelian fourfolds we prove:

\begin{enumerate} 

\item If $(A,L)$ is generic of type $(d_1,d_2,d_3,d_4)$ with $h^0(L)\geq 31$ and $$
\begin{array}{rl}
{(d_1,d_2,d_3,d_4)\not\in }&{\{(1,1,2,16),(1,2,2,8), (1,2,4,4),}\\ 
{}&{\,(1,3,3,6),(2,2,2,4), (2,2,4,4), (2,4,4,4)\}}\\
\end{array}
$$ then $L$ is normally generated.

\item If $L$ is of type $(1,1,2,16)$, $(1,2,2,8)$, $(1,2,4,4)$, $(2,2,2,4)$, $(2,2,4,4)$ or $(2,4,4,4)$ then it is never normally generated.

\end{enumerate}

Now the open case is when $L$ is of type $(1,3,3,6)$.

{\bf Acknowledgements: } I thank K. Hulek for introducing me in these topics.

{\bf Notations:} We will work over the field of the complex numbers. We will follow the definitions and notation of the book \cite{bila}.

\section{The main Theorem.}

Let $(B,M_0)$ be a principal polarized variety with $B=\CC^g/\Lambda$ and $M_0=L(H,\chi_0)$. Suppose that $M_0$ is of characteristic $0$. Let $\lambda_1,\ldots,\lambda_g,\mu_1,\ldots,\mu_g$ be the corresponding  symplectic base. 

We will denote by $V_1$ and $V_2$ the real spaces generated respectively by $\lambda_1,\ldots,\lambda_g$ and $\mu_1,\ldots,\mu_g$. Any $x\in \CC^g$ decomposes in a unique way  $x=x_1+x_2$ with $x_i\in V_i$.  Let $\Lambda_i=V_i\cap \Lambda$.

Let $D=(d_1,\ldots,d_g)$ with $d_i|d_{i+1}$ be a type. Let $$\Lambda'=\langle \lambda_1,\ldots,\lambda_g,d_1\mu_1,\ldots,d_g\mu_g \rangle.$$ We consider the abelian variety $A=\CC^g/\Lambda'$. There is a natural projection $p:A\lrw B$. It is clear that the line bundle 
$L=p^*M_0$ on $Y$ has type $D$ and  characteristic $0$.

Let $K_1=\{\frac{n_1}{d_1}\lambda_1+\ldots+\frac{n_g}{d_g}\lambda_g|\,n_i\in \ZZ\, ,1\leq i\leq g\}/\Lambda_1$. We have the following isomorphisms (see \cite{iy}):
$$
H^0(L)=\sum_{w\in K_1}H^0(t^*_w M_0); \qquad
H^0(L^2)=\sum_{w\in K_1}H^0(t^*_{\frac{w}{2}} M_0^2).
$$

\begin{teo}\label{iyer1}

The map $H^0(L)\otimes H^0(L)\lrw H^0(L^2)$ factorizes through the maps:
$$
\sum_{u\in K_1}H^0(t_{u}^*M_0)\otimes H^0(t_{-u+w}^*M_0)\lrw H^0(t_{\frac{w}{2}}^*M_0^2)
$$
where $w$ runs over $K_1$.

\end{teo}
{\bf Proof:} See \cite{iy}. \qed 

Let us study explicitly the surjectivity of these factorizations. We will denote by $\xi^c_w$ the canonical theta functions for the line bundle $M_0$ and by $\vartheta^c_w$ the canonical theta functions of $M_0^2$. Let $u,w\in K_1$. We know that  $H^0(t_{u}^*M_0)$ and $H^0(t_{-u+w}^*M_0)$ are generated respectively by the functions $\xi_u^0$ and $\xi_{-u+w}^0$. Thus, from the Multiplication Formula (\cite{bila},7.1.3):
$$
\xi_u^0 \xi_{-u+w}^0=\sum_{z\in Z_2} \vartheta_{u+z-\frac{w}{2}}^0(0) \vartheta^{\frac{w}{2}}_z
$$
where $Z_2=B_2\cap V_1$. Now, we can obtain the following result:

\begin{teo}\label{principal}

The line bundle $L$ is $2$-normal if and only if for each $w\in K_1/2K_1$ the matrix:
$$
\left(\vartheta_{u+z-\frac{w}{2}}^0(0)\right)_{u\in K_1,z\in Z_2}
$$
has rank $2^g$.

\end{teo}
{\bf Proof:} From the previous arguments and Theorem \ref{iyer1} it follows that $L$ is $2$-normal exactly when the matrices 
$$
\left(\vartheta_{u+z-\frac{w}{2}}^0(0)\right)_{u\in K_1,z\in Z_2}
$$
have rank $2^g$ for any $w\in K_1$. But, it is clear that if we take $w,w'\in K_1$ such that $w-w'\in 2K_1$ then the two corresponding matrices have the same rank. \qed

\section{Cases where the projective normality fails.}

In this section we will use the main Theorem to give two situations where the projective normality fails. We will suppose that there is a $d_j$ equal to $2$. In particular, let $d_{i_0}=2$ where $i_0$ is the minimum index verifying this condition. Take $w=\lambda_{i_0}/2$. We decompose $K_1$ in two subgroups:
$$
\begin{array}{l}
{K_{11}=\{\frac{n}{2}\lambda_{i_0}|\, n\in \ZZ\}/\Lambda_1}\\
{K_{12}=\{\frac{n_{i_0+1}}{d_{i_0+1}}\lambda_1+\ldots+\frac{n_g}{d_g}\lambda_g|\,n_i\in \ZZ\, ,i_0+1\leq i\leq g\}/\Lambda_1}\\
\end{array}
$$
We have the following lemma:

\begin{lemma}
With the previous notation, the  matrices:
$$
\left(\vartheta_{u+z-\frac{w}{2}}^0(0)\right)_{u\in K_1,z\in Z_2}\hbox { and }
\left(\vartheta_{u+z-\frac{w}{2}}^0(0)\right)_{u\in K_{12},z\in Z_2}
$$
have the same rank.
\end{lemma}
{\bf Proof:} Let $u=u_1+u_2$, with $u_1\in K_{11}$ and $u_2\in K_{12}$. Suppose that $u_1=\frac{\lambda_{i_0}}{2}$. Then, since $M_0$ is of characteristic $0$ and by using the properties of the canonical theta functions: 
$$
\vartheta_{u+z-\frac{w}{2}}^0(0)=\vartheta_{\frac{\lambda_{i_0}}{2}+u_2+z-\frac{\lambda_{i_0}}{4}}^0(0)=\vartheta_ {-\frac{\lambda_{i_0}}{2}-u_2-z+\frac{\lambda_{i_0}}{4}}^0(0)=
\vartheta_{-u_2+z-\frac{w}{2}}^0(0)
$$ \qed

By this result we will center our attention in the matrix:
$$
Q=\left(\vartheta_{u+z-\frac{w}{2}}^0(0)\right)_{u\in K_{12},z\in Z_2}
$$
with $w=\frac{\lambda_{i_0}}{2}$.

We decompose the group $Z_2$ in :
$$
\begin{array}{l}
{Z_{21}=Z_2\cap \langle \frac{\lambda_{i_0}}{2} \rangle}\\
{Z_{22}=Z_2\cap \langle \frac{\lambda_{1}}{2},\ldots,\frac{\lambda_{i_0-1}}{2},\frac{\lambda_{i_0+1}}{2},\ldots,\frac{\lambda_{g}}{2} \rangle}\\
\end{array}
$$
this induces a decomposition for the matrix $Q$:
$$
Q=\left(\matrix{Q_1 & Q_2}\right)
$$
with
$$
Q_1=\left(\vartheta_{u+z-\frac{w}{2}}^0(0)\right)_{u\in K_{12},z\in Z_{22}},\qquad Q_2=\left(\vartheta_{u+z+\frac{\lambda_{i_0}}{2}-\frac{w}{2}}^0(0)\right)_{u\in K_{12},z\in Z_{22}}
$$

Let $u\in K_{12}, z\in Z_{22}$. Then:
$$
\begin{array}{rl}
{\vartheta_{u+z-\frac{w}{2}}^0(0)}&{=\vartheta_{u+z-\frac{\lambda_{i_0}}{4}}^0(0)=\vartheta_{-u-z+\frac{\lambda_{i_0}}{4}}^0(0)=}\\
{}&{=\vartheta_ {u+(-2u-z)+\frac{\lambda_{i_0}}{2}-\frac{\lambda_{i_0}}{4}}^0(0)}\\
\end{array}
$$
If $u$ is a $4$-torsion point of $K_{12}$ then $-2u-z$ is an element of $Z_{22}$. Therefore, in this case, the sum of the columns of $Q_1$ is the same that the sum of the columns of $Q_2$. We deduce that the matrix $Q$ has not maximal rank $2^g$. Thus, we have prove the following Theorem:

\begin{teo}\label{fail1}

Let $(A,L)$ be a $g$-dimensional polarized abelian variety of type $(d_1,\ldots,d_g)$, with $d_j\leq 4$, $1\leq j\leq g$ and some $d_i=2$. Then $L$ is not normally generated. \qed

\end{teo}

On the other hand, we consider the matrix:
$$
Q_1-Q_2=\left(\vartheta_{u+z-\frac{w}{2}}^0(0)-\vartheta_ {u+z+\frac{\lambda_{i_0}}{2}-\frac{w}{2}}^0(0)\right)_{u\in K_{12},z\in Z_{22}}
$$
We decompose $K_{12}$ in the following way:
$$
K_{12}=K_{12}^0\cup K_{12}^1
$$
where $K_{12}^0$ are the $2$-torsion points of $K_{12}$ and $K_{12}^1$ the rest of them. Each  $u\in K_{12}^1$ has a unique inverse $-u$ on $K_{12}^1$.  Furthermore:

- If $u\in K_{12}^0$ then $(Q_1-Q_2)_{u,z}=0$ for any $z\in Z_{22}$, because:
$$
\vartheta_{u+z+\frac{\lambda_{i_0}}{2}-\frac{\lambda_{i_0}}{4}}^0(0)=\vartheta_{-u-z-\frac{\lambda_{i_0}}{2}+\frac{\lambda_{i_0}}{4}}^0(0)=\vartheta_{u+z-\frac{\lambda_{i_0}}{4}}^0(0)
$$

- If $u\in K_{12}^1$ then $(Q_1-Q_2)_{u,z}=-(Q_1-Q_2)_{-u,z}$ for any $z\in Z_{22}$, because:
$$
\begin{array}{rl}
{\vartheta_{u+z-\frac{w}{2}}^0(0)-\vartheta_{u+z+\frac{\lambda_{i_0}}{2}-\frac{w}{2}}^0(0)}&{=\vartheta_{-u-z+\frac{w}{2}}^0(0)-\vartheta_{-u-z-\frac{\lambda_{i_0}}{2}+\frac{w}{2}}^0(0)=}\\
{}&{=\vartheta_{-u+z+\frac{\lambda_{i_0}}{2}-\frac{w}{2}}^0(0)-\vartheta_{-u+z-\frac{w}{2}}^0(0)}\\
\end{array}
$$

We see that   $Q_1-Q_2$ has at least $\#K_{12}^0$ null rows and $\frac{1}{2}\#K_{12}^1$ pairs of rows with opposite sign. Let $h^0(L)=2^{g+1-{i_0}} n$, then $\#K_{12}=2^{g-{i_0}} n$ and $\#K_{12}^0=2^{g-{i_0}}$. From this:
$$
rank(Q_1-Q_2)\leq \frac{1}{2}\#K_{12}^1=\frac{1}{2}(\#K_{12}-\#K_{12}^0)=2^{g-{i_0}-1} n-2^{g-{i_0}-1}
$$
When this rank is less than $2^{g-1}$ then $Q$ has not maximal rank, so $L$ is not normally generated:
$$
rank(Q_1-Q_2)<2^{g-1}\iff n<2^{i_0}+1 \iff n\leq 2^{i_0}\iff h^0(L)\leq 2^{g+1}
$$
Since a necessary condition for the projective normality is $h^0(L)\geq 2^{g+1}-1$, we have proof the following heorem:

\begin{teo}\label{fail2}

Let $(A,L)$ be a $g$-dimensional polarized abelian variety of type $(d_1,\ldots,d_g)$, with some $d_i=2$. If $h^0(L)=2^{g+1}$, then $L$ is not normally generated. \qed

\end{teo}

\section{Projective normality for $g=3$ and $g=4$.}

In this section we will use the main Theorem \ref{principal} to study the projective normality of abelian threefolds and fourfolds making explicit computations. Since this is an open condition we will give explicit examples where the projective normality holds. This will prove the projective normality in the generic case.

To apply the Theorem \ref{principal}, we have to compute $\vartheta_{c}^0(0)$ with $c\in V_1$. We will write this functions by using the classical theta functions. We will follow the Chapter 8 of \cite{bila}. Let $Z\in M_g(\CC)$ be a symmetric matrix with $Im(Z)>0$ and let $c\in \CC^g$ such that $c=Zc^1$, with $c^1\in \RR^g$. By the Lemma (8.5.2., \cite{bila}) and the properties of the theta functions we have: 
$$
\begin{array}{rl}
{\vartheta_{c}^0(0)}&{=\vartheta_{0}^c(0)=\vartheta\left[\matrix{c^1 \cr 0}\right](0,Z)=}\\
{}&{=\e{\pi \imath ^tc^1Zc^1}\vartheta\left[\matrix{0 \cr 0}\right](Zc^1,Z)=\e{\pi \imath ^tc^1Zc^1}\theta(Z,Zc^1)}\\
\end{array}
$$
where $\theta$ is the Siegel theta function:
$$
\theta(Z,v)=\sum_t \e{\pi \imath ^tt Z t +2 \pi\imath ^t t v}
$$
Now, we can rewrite the Theorem \ref{principal}.
Let $I=\{(\frac{n_1}{d_1},\ldots,\frac{n_g}{d_g})|\,0\leq n_i< d_i\, , 1\leq i\leq g\}$, $I'=\{(\frac{n_1}{d_1},\ldots,\frac{n_g}{d_g})|\,0\leq n_i< m.c.d.(d_i,2)\, , 1\leq i\leq g\}$ and $J=\{(\frac{a_1}{2},\ldots,\frac{a_g}{2})|,0\leq a_1,\ldots,a_g\leq 1\}$. Then:

\begin{teo}

The line bundle $L$ is $2$-normal if and only if for any $w\in I'$ the matrices:
$$
\left(\e{\pi\imath^t (i+j-\frac{w}{2})Z(i+j-\frac{w}{2})}\theta(Z,Z(i+j-\frac{w}{2}))\right)_{i\in I,j\in J}
$$
have rank $2^g$. \qed

\end{teo}

The classical Theta functions can be evaluated with the help of a computer. In particular Mathematica 5.0 (\cite{wo}) includes an implementation of the Siegel theta function. Furthermore, in \cite{deheboho} algorithms for their computation  are given. However, one can simplify the situation working on a particular case. We can suppose that $Z=X+Y$ where $X$ is a symmetric integer matrix and $Y=k Id$, with $k\in \CC$ and $Im(k)>0$. Now, by the Theta Transformation Formula (8.6.1, \cite{bila}):
$$
\e{\pi \imath ^tc^1Zc^1}\vartheta\left[\matrix{0 \cr 0}\right](Zc^1,Z)=\e{\pi \imath ^tc^1Zc^1}\vartheta\left[\matrix{0 \cr 0}\right](Zc^1,Y)
$$
Because $Y$ is a diagonal matrix, the functions $\vartheta\left[\matrix{0 \cr 0}\right](Zc^1,Z)$ can be decomposed in a product of elliptic functions. 

We have applied this method for $g=3$, and the decomposition:
$$
X=\left(\matrix{0 & 0 & 1 \cr 0 & 0 & 2 \cr 1 & 2 & 0}\right);\qquad
Y=(1+\sqrt{1/3} \imath) Id $$
We have checked all the possibilities for the type $(d_1,d_2,d_3)$ when $7\leq d_1d_2d_3\leq 2^3\cdot 3!$. In other case, the results of Iyer \cite{iy} provide the projective normality. We have obtained the following result:

\begin{teo}

Let  $(A,L)$ be a generic abelian threefold of type $(d_1,d_2,d_3)$ with $h^0(L)\geq 15$. If  $
(d_1,d_2,d_3)\not\in \{(1,2,8),(1,3,6),(2,2,4),(2,4,4)\}$ then $L$ is normally generated.

\end{teo}

Moreover, applying the Theorems \ref{fail1} and \ref{fail2} for $g=3$ we obtain the following corollary:

\begin{cor}

An abelian threefold $(A,L)$ of type $(1,2,8),(2,2,4)$ or $(2,4,4)$ is never projectively normal. \qed
 
\end{cor}

For $g=4$ we have used the  decomposition:
$$
X=\left(\matrix{0 & 0 & 0 & 1 \cr 0 & 0 & 0 & 2 \cr 0 & 0 & 0 & 3 \cr 1 & 2 &3 & 0}\right);\qquad
Y=(1+\sqrt{1/3} \imath) Id
$$
Now, we have checked all the possibilities for the type $(d_1,d_2,d_3,d_4)$ when $15\leq d_1d_2d_3d_4\leq 2^4\cdot 4!$. We have obtained the following results:

\begin{teo}
Let  $(A,L)$ be a generic abelian fourfold of type $(d_1,d_2,d_3,d_4)$ with $h^0(L)\geq 31$. If  $$
\begin{array}{rl}
{(d_1,d_2,d_3,d_4)\not\in }&{\{(1,1,2,16),(1,2,2,8), (1,2,4,4),}\\ 
{}&{\,(1,3,3,6),(2,2,2,4), (2,2,4,4), (2,4,4,4)\}}\\
\end{array}
$$ then $L$ is normally generated. \qed
\end{teo}

In this case, from  Theorems \ref{fail1} and \ref{fail2} it follows:

\begin{cor}

If $(A,L)$ is an abelian fourfold of type $(1,1,2,16)$, $(1,2,2,8)$, $(2,2,2,4)$, $(2,2,4,4)$ or $(2,4,4,4)$ then $L$ is not normally generated. \qed

\end{cor}

Note, that the open cases for $g=3$ and $g=4$ are respectively the types $(1,3,6)$ and $(1,3,3,6)$. We have not found examples of projectively normal polarized abelian varieties of theses types. On the other hand, we have explored some other cases for $g>4$. Thus, we make the following conjectures:

\begin{conj}

Let $(A,L)$ be a polarized abelian variety of dimension $g$. If $L$ is of type $(1,3,\cdots,3,6)$ then $L$ is not normally generated. \qed

\end{conj}

\begin{conj}

Let $(A,L)$ be a polarized abelian variety of dimension $g$. If $L$ is of type $(1,\cdots,1,d)$ with $d\geq 2^{g+1}-1$ then $L$ is  normally generated. \qed

\end{conj}

\end{document}